\documentclass[12pt]{article}
\usepackage{theorem}
\theorembodyfont{\rmfamily}
\usepackage[dvips]{graphics}
\newcommand{\fig}[3]{
\begin{figure}
\center{\includegraphics{#1}}
\caption{#3}
\label{#2}
\end{figure}
}
\newtheorem{thm}{Theorem}[section]
\newtheorem{lemma}[thm]{Lemma}
\newtheorem{defn}[thm]{Definition}
\newtheorem{formula}[thm]{Formula}
\newtheorem{remark}[thm]{Remark}
\newtheorem{question}[thm]{Question}
\newtheorem{cor}[thm]{Corollary}
\newcommand{\Bbb}[1]{%
{\bf #1}}
\newcommand{\frak}[1]{%
{\bf #1}}
\newcommand{\qed}{{\bf q.e.d}}

\title{Random Delaunay triangulations \\ and metric uniformization}
 \author{Gregory Leibon }

\begin{document}

\maketitle

\begin{section}{Introduction}

The primary goal of this paper is to develop a new connection
between the discrete conformal geometry problem of disk pattern
construction and the continuous conformal geometry problem of metric
uniformization. 
In a nutshell, we discuss how to construct disk  patterns 
by optimizing an
objective function, which turns out to be intimately related to
hyperbolic volume. 
With the use of random Delaunay triangulations we 
then average this  objective function  to construct an 
objective function on the  metrics 
conformal to a fixed
one. Finally using this averaged objective function we 
may reprove the uniformization theorem in two dimensions. 

This paper is organized as follows.  In section \ref{sec:2} we 
introduce 
the use of random Delaunay triangulations by presenting a new 
proof of the Gauss-Bonnet formula. 
In section \ref{sec:3} we present the
disk pattern ideas and the objective function in the disk pattern  
setting.
In section \ref{sec:4} we use the techniques developed in section
\ref{sec:2}  to average the objective function from section \ref{sec:3}
and  produce the objective function on metrics.
In section \ref{sec:5} we point to some open questions.

This article is based on the author's thesis \cite{Le},
where readers can find
a  detailed treatment of everything that takes place here.

\end{section}

\begin{section}{A  Random Proof of the Gauss-Bonnet Theorem}\label{sec:2}

Recall that the  Gauss-Bonnet  formula  may be stated 
for a compact  boundaryless Riemannian surface $M$  as: 
\[ \frac{1}{2\pi} \int_{M} k dA = \chi (M) \]
with $k$ the Gaussian curvature and $\chi(M)$ the surface's 
Euler characteristic. The Euler
characteristic is a topological invariant 
which can be computed relative to any triangulation  via 
\[ \chi(M) = F-E+V\]
with $F$, $E$, and $V$ the number of faces, edges and vertices in the
triangulation. 
The proof here is accomplished by
randomly triangulating the surface and then noting that while 
$\chi(M)$   is constant $F$, $E$, and $V$ 
are  now random variables and  have expected values which can be computed. As the density of the randomly distributed 
vertices goes to infinity  one finds these expected values produce the 
 Gauss-Bonnet formula, along with a probabilistic interpretation 
of Gaussian curvature.

To begin with, we need to define what we mean by a random 
triangulation of a surface. The first step is to ignore the fact that 
anything random  is going on here and to simply attempt to construct  
a geodesic triangulation in a fixed metric $g$
from a given set of points,
${\bf p} = \{ p_1, ... ,p_n \}$. To accomplish this we 
produce 
an abstract two complex 
by examining all the triples  and pairs in $\{ p_1, ... ,p_n \}$ 
and deciding whether or not to put in a face for a given triple or 
an edge for a given pair.
 This decision procedure 
will be relative to a 
certain positive number $\delta$ -- the {\it  decision radius}.  The procedure 
is to put in a face for a triple or an edge for a pair if the 
triple or pair lies on a disk of 
radius $<  \delta$ which has its interior empty of points 
in $\{ p_1, ... ,p_n \}$.   If a triple is on an empty disk then each
pair in this triple is on an empty disk, so we indeed have a
2-complex.  
 This proceedure 
is called Delaunay's ``empty sphere'' method  
and was introdueced in \cite{Da}.
It is elementary to see that there is a positive decision radius
such that  one can view the edges of this 2-complex 
as geodesics in the $g$ metric. 
When this  procedure forms a triangulation of $M$  we call  the 
resulting triangulation the {\it Delaunay triangulation}.   

At this point it is useful to introduce a geometric criterion on a 
set of points
 ${\{ p_1, ... ,p_n \}}$ guaranteeing that it produces a 
Delaunay triangulation.

\begin{defn}
We  will call a set of points  ${\{ p_1, ... ,p_n \}}$ {\it generically 
$\delta$-dense} 
if each open ball of radius $\delta$ contains at least one $p_i$ and 
${\{ p_1, ... ,p_n \}}$ contains no four    points
on a disk of radius less than
$\delta$. 
\end{defn}

It is straightforward to see
that
  
\begin{lemma}
\label{dtri}
There is a $\delta>0$ such that if ${\{ p_1, ... ,p_n \}}$ is
generically $\delta$-dense,
then ${\{ p_1, ... ,p_n \}}$ forms a Delaunay triangulation. 
\end{lemma}

Now enters the randomness.  From lemma \ref{dtri}, when our
points are distributed with a high density we expect they will typically
form Delaunay triangulations. This procedure is somewhat independent
of our choice of distribution, but  for all that takes place here we
will assume that  we are using a Poisson point process 
relative to a density denoted $\lambda$ (see remark
\ref{unifrem}).  Letting  $ \Bbb E_{\lambda}(L)$   
denote  the expected value of a random variable $L$  and letting
$O(\lambda^{-\infty})$  mean a quantity decaying faster than any 
polynomial in $\lambda$,
we indeed have that...
\begin{lemma}\label{size} 
The probability that a set of points form a triangulation is $1+
O(\lambda^{-\infty})$.    Also, if $L$ is any one of the random variables 
$E$, $V$, or $F$, then
$L$  has expected value 
equal to  $ \Bbb E_{\lambda} \left(\overline{L} \right)
 + O(\lambda^{-\infty}) $,  with $\overline{L}$ the random variable which is $L$
when the points form a triangulation and zero otherwise.   
\end{lemma}

From this lemma
 one sees that the constant $\chi(M)$ satisfies
\[\chi (M) =  \Bbb E_{\lambda}\left(\overline{\chi(M)} \right) +   O(\lambda^{-\infty}) = 
\Bbb E_\lambda \left( \overline{V} - \overline{E} + \overline{F}\right)  + 
O(\lambda^{-\infty}).
\]  
Applying the
fact  that expected 
values add gives us
\[ \chi (M) = \Bbb E_\lambda \left(\overline{V}\right) -   \Bbb E_\lambda
\left(\overline{E}\right)  + 
\Bbb E_\lambda
\left(\overline{F}\right)  +  O(\lambda^{-\infty})  . \]
Furthermore, in an actual  triangulation we have $\frac{3}{2}E = F$ so
\[\chi (M) =  \Bbb E_\lambda \left( \overline{V} \right) - \frac{1}{2} \Bbb E_\lambda
\left(\overline{F}\right)  
 +  O(\lambda^{-\infty}). \]
Since $\lambda$ is the density, if we let $A$ denote the area of $M$ then the expected
number of vertices is
$\lambda A
$.  This observation along with
lemma $\ref{size}$ gives us 
\[\chi (M)  = A \lambda  - \frac{1}{2}  \Bbb E_\lambda (F) 
 +  O(\lambda^{-\infty}) . \]
So we have... 
\begin{formula}[Euler-Delaunay-Poisson Formula]
\label{edpf}
\[  \chi (M) = \lim_{\lambda \rightarrow \infty} 
\left( A \lambda  - \frac{1}{2}  \Bbb E_\lambda (F)  \right). \]
\end{formula}

The goal now becomes to compute $\Bbb E_\lambda (F)$. 
Let $\delta$ be small enough to satisfy lemma \ref{dtri}, let   $V_{\delta}
\subset  \times^3 M=  M \times M \times M $   be the set of ordered 
triples living on circles of radius
less than $\delta$, and let $a(y)$ be the area of the disk associated to the
triple  $y \in V_{\delta}$. If you are familiar with how to 
compute using  the Poisson process you will find
\begin{eqnarray} 
\label{ex1}
 \Bbb E_{\lambda}(F) =   
 \frac{1}{6} \int_{V_{\delta}} e^{-\lambda a(y)} (\lambda dA)^3. 
\end{eqnarray}
(If you are not familiar with this situation you should see remark \ref{sec}.)

In order to explicitly compute the integral in formula \ref{ex1} 
it is necessary to put
coordinates  on $V_{\delta}$.  To accomplish this, first one chooses 
a way to discuss
directions  at all but a finite number of tangent planes of $M$ (via
an orthonormal frame field).
Then one can parameterize a full measure subset of 
$V_{\delta}$ with a  subset
of  
$\{ (\theta_1, \theta_2,\theta_3, r,p) \in S^1 \times S^1 \times  S^1 
\times (0 , \delta) \times M \}$ 
by starting  at the point $p \in M$  
and moving a distance $r$ in each of the three directions $ \theta_1$, 
$\theta_2$, and $\theta_3$.
Note that when fixing $p$ and $r$ and varying 
$\theta_i$  we 
produce a Jacobi field, whose norm we shall denote   $j_i$.
Using this notation,
letting   $d \vec{\theta} =  d \theta_1  d \theta_2 d \theta_3$,
and 
letting  $ \nu( \vec{\theta})$ be the area of the triangle 
in the Euclidean unit circle 
with vertices at the points corresponding to the $\{ \theta_i\}$, 
a straightforward computation  shows us that formula \ref{ex1} is   
\begin{eqnarray}
\label{ex2}
 \frac{\lambda^3}{6}  \int_{M} \int_0^{\delta} \int_{\times^3 S^1} 
e^{-\lambda a(\vec{\theta},r,p)} j_{\theta_1}j_{\theta_2}j_{\theta_3}  
\nu( \vec{\theta}) d
\vec{\theta} dr dA.
\end{eqnarray}
in these coordinates.

The Taylor expansions of  $a(\vec{\theta},r,p)$ and
$j_{\theta}$ are controlled by the Gaussian 
curvature up to the fourth and third order
terms respectively.  So, after potentially shrinking $\delta$ a bit to 
 exploit this
control, we may Taylor expand, integrate, and apply the mean value
theorem to express
\ref{ex2} as
\begin{eqnarray}
\label{ex3}
  \Bbb E_{\lambda} (F) =  2 A \lambda -    
\frac{1}{ \pi } \int_{M} k dA  + O\left(\lambda^{-\frac{1}{2}}\right). 
\end{eqnarray}

In particular, formula \ref{edpf} may now be plugged
into the Euler-Delaunay-Poisson Formula 
 to give simultaneously  a probabilistic 
interpretation of curvature and a proof of the Gauss-Bonnet theorem.
\begin{thm}[Euler-Gauss-Bonnet-Delaunay Formula]
\label{egbdf}
\[\chi (M) =  \lim_{\lambda \rightarrow \infty} \left( A \lambda  - \frac{1}{2}  \Bbb
E_{\lambda} (F)\right) = \frac{1}{2 \pi} \int_{M} k dA .\]
\end{thm}

Note this allows us to interpret the Gaussian curvature as the 
density of the 
defect in the expected number of faces in a random Delaunay triangulation
in the surface's geometry 
relative to what would be expected in 
the Euclidean plane. 

\begin{remark}
It should be noted that the  lemmas and computations above are completely
elementary with the exception of the use of the following bit of
geometry.

\begin{lemma}[The Small Circle Intersection Lemma]
\label{lit}
There is a $\delta >0$ such that if a triple of points lies on 
 the boundary of a
disk with  radius less than $\delta$, then this disk is unique among disks of 
radius less than  $\delta$. 
\end{lemma}
This lemma is at the core of the proof of lemma \ref{dtri}, as well 
the justification for the well-definedness of $a(y)$ and the 
parameterization in formula \ref{ex2}.   A couple of proofs of this lemma can 
be found in
\cite{Le}, where it is shown that we can get an explicit grip on the
necessary
$\delta$ by using any
$\delta
< min\{\frac{i}{6},\tau\}$, where $i$ is the surface's injectivity radius and
$\tau$ the surface's strong convexity radius.
\end{remark}

\begin{remark}\label{unifrem} 
Any  reasonable form of point distribution will fare as
well as the Poisson point process with regards to this entire proof.  
This includes the uniform 
distribution of $n$ points, which is simply the 
product measure on
$M \times \ldots \times M = \times^nM$.    
This uniform  choice in fact eliminates all but the most basic
probabilistic thinking, though not without a certain aesthetic sacrifice.
\end{remark}

\begin{remark}\label{sec}
Here we will give an informal  description of the Poisson point
process and  the 
derivation of formula \ref{ex1}.  To describe the Poisson point
process, imagine breaking the surface up 
into pieces of size
$dA$ small enough so that the probability 
that one of these pieces contains more
than one point is negligible.  Denote one of these  little  chunks by 
$q$ and let 
$X_q$ be 1 if the region $q$ contains a point and 0 otherwise.  The heart
of the Poisson point process is the assumption  that the
$X_q$ are independent and  the probability that $q$ contains a point
is $\Bbb{E}_{\lambda} (X_q) = \lambda dA$. 
Exploiting this  independence we see that
the probability of a region  $R$ of area
$A$ being empty of points is 
\[ \Bbb{E}_{\lambda}(\Pi_{q \in R} X_q)  \approx   
(1-\lambda dA)^{\frac{A}{dA}} \approx e^{-\lambda A}.\]
 Let $R_t$ be the function which is one if the ``disk''  formed by a ``triple''
 $t = \{q_1,q_2,q_3\}$  is empty of other points, 
and zero otherwise (the quotations are used because it is only after
taking the limit that we really have a triple of points and an
actual disk). In
particular note that $R_t$ has expected value 
$e^{- \lambda a(t)}$.   Hence by independence we have that
 \[ \Bbb E_\lambda (R_{\{q_1,q_2,q_3\}} X_{q_1} X_{q_2} X_{q_3}) =
 e^{-\lambda a(t)} (\lambda dA)^3 \] 
is   the probability of a ``triple'' forming
a face. 
In particular 
 by the linearity of expected values
\[\Bbb E_\lambda (F) =  \sum_{``triples"}  \Bbb E_\lambda (R_t) = 
 \sum_{``triples"} e^{a(t)} (\lambda dA)^3,\] 
which gives us  formula \ref{ex1} in the limit. 
\end{remark}

\end{section}

\begin{section}{Discrete Conformal Uniformization}\label{sec:3}

The goal here is to describe the needed disk pattern ideas.
In order to facillitate the use  of these ideas in the metric
world we will
present our discrete objects as natural data living within  actual
geodesic triangulations of Riemannian surfaces.
For example, imagine starting with 
a triangulation of a hyperbolic surface, by which we  mean a
geometric surface 
having constant Gaussian curvature $-1$.  
 Let
 $s$ and $t$ be two triangles  in  this  triangulation 
which share an  edge $e$, and denote  the complement of the intersection 
angle  between  the disks in
 which  $s$ and $t$  are inscribed  as $\psi^{e}$. See Figure 1. 
The key is to note that 
$\psi^e$ can be written down 
in terms of the angles within the triangles as

\fig{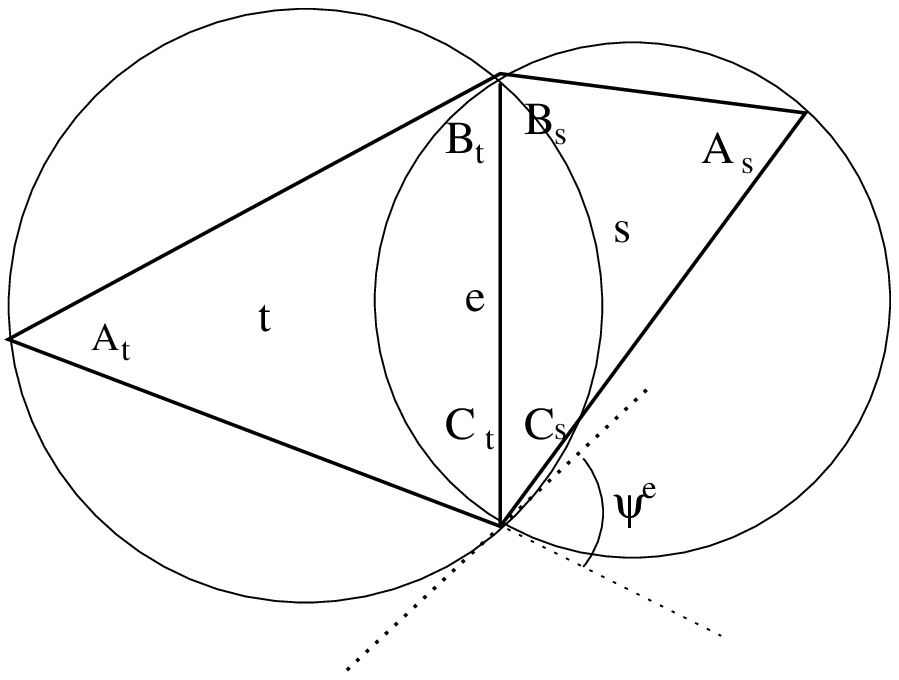}{nota1}{The notation.}

\begin{formula} \label{theta} 
    \[\psi^e= \psi^e_{s} + \psi^e_{t}  \]
with 
\[ \psi^e_{\tau} =  \frac{B_\tau + C_\tau - A_\tau}{2}\]
using the notation in  the first figure.  
\end{formula}

{\bf Comment on Proof:} 
This formula always holds on a constant curvature surface, and
a  proof  can be  found in \cite{Les}.  For our purposes
 it is  useful  to introduce the key object showing up in the
 proof of  the negative curvature case. 
This object is an
ideal hyperbolic prism.   The prism  is constructed from the angle data of
a  hyperbolic triangle, namely a set of positive angles 
$\{A,B,C\}$ such that $A + B + C -\pi <0$.
To construct it first form a hyperbolic 
triangle with the $\{A, B, C\}$ data,  then  place the triangle on
 a hyperbolic plane and 
      then place the plane in hyperbolic three space.  Now 
union this triangle with   the geodesics 
       perpendicular to this 2-plane going through the 
vertices of the triangle.  The prism of interest is
 the convex hull of this 
arrangement. See Figure 2.

\fig{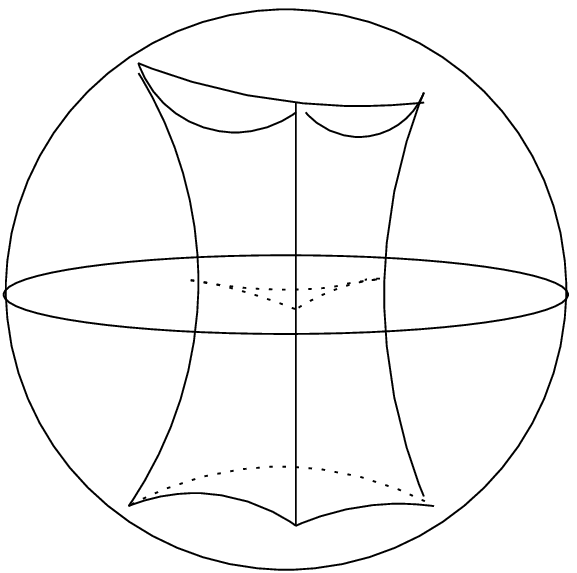}{prism}{The prism.}

Now suppose we have any  triangulation of a Riemannian surface with
varying curvature;
What  data may we  collect?  Well, we may collect the topological
triangulation  $\frak{T}$ and all the triangle angles.
From formula \ref{theta} this triangle angle data may be organized into 
the collection of partial angles $\psi^e_t$,  whose 
values we will identify with the coordinates in a
$3 F$ dimensional real vector. 
It should be pointed out that  
the angles inside the
triangles are linearly determined by the  $\psi^e_t$ via inverting 
formula \ref{theta}.   To be the data of a  Delaunay
triangulation
forces certain natural linear constraints, which we
will capture with the following definition.

\begin{defn}
Let an {\it angle sytem} be any point $x$ in the above vector space
such that
$x$'s associated 
triangle angles  are in $ (0,\pi)$ and  sum up to $2 \pi$ at any
vertex.  A {\it Delaunay angle system} is an angle system $x$ where
the complement of the 
{\it informal intersection angle}, 
$\psi^e(x) =  \psi^e_{s}(x) + \psi^e_{t}(x)$, 
has values in $(0,\pi)$.
\end{defn}

We call two  angle systems $x$ and $y$ {\it conformally equivalent}
if they  share 
 the same  informal intersection angles.
We record this with idea the following definition. 

\begin{defn}
$x$ and $y$ will be called  conformally equivalent if 
$\psi^e (x) = \psi^e(y)$
for all
$e$.
\end{defn}

Motivated by the  Gauss-Bonnet formula we also define the following.

\begin{defn}
The {\it curvature} of a triangle $t$ relative to an angle system 
$x$ is defined to be $\pi$
subtracted from sum of the triangle angles in $t$ determined by $x$.
If every 
triangle has negative
curvature relative to $x$, then $x$ will be said to have {\it negative
curvature}.
The set of Delaunay angle systems
with negative curvature will be denoted $\frak{N}_{x}$.
\end{defn}

There  is a simple set of linear equations 
which will gaurentee that a 
Delaunay angle system is conformally equivalent 
an angle system in   $\frak{N}_{x}$.
Let $S$ be a set of triangles in $\frak{T}$, and denote the 
cardinality of $S$ as $|S|$. 

\begin{defn}
An angle system $x$ will be called {\it teleportable} if 
for any set of 
triangles $S$ we have 
\[  \sum_{e \in S} (\pi - \psi^e(x)) > \pi |S|. \]
\end{defn}

As an immediate  consequence of 
theorem 2  in  \cite{Les}
we have the following lemma.

\begin{lemma}[The Discrete Teleportation Lemma]
\label{dtele}
A Delaunay angle system $x$ is teleportable  
if and only if $x$ is conformally 
equivalent to a negative curvature  Delaunay angle system.
\end{lemma}

The proof of this lemma  is completely linear in nature. 

In the
setting here we will only be concerned with Delaunay  angle systems 
satisfying  the
conditions of this lemma, which  includes the angle data associated to a 
Delaunay triangulation of surface with varying negative curvature.

The  pleasure derived from the Delaunay angle systems comes from  
 a wondrous objective function which lives on  $\frak{N}_{x}$.  Let $V_t(x)$ 
denote the volume of the ideal hyperbolic 
prism constructed from  a triangle $t$'s angle data  relative to $x$,
as described in the previous section.
Now simply let the objective function  be
\[H(x) = \sum_{t \in \frak{T}} V_t(x). \]
The wonder of this function  can  best be felt by examining its
differential.  To compute
this    use  $x$'s
angle data in $t$ to construct a hyperbolic triangle
and let $l^e_t(x)$ denote the length of the edge $e$ in this 
triangle. In \cite{Les}  
the following formula is produced:

\begin{formula} \label{dif}
\[ dH = - \sum_{(e,t)}  \log\left(\frac{\cosh(l^e_t(x))-1}{2}\right) 
d \psi_t^e\]
\end{formula}

The  tangent space at any point of
 $\frak{N}_{x}$ is precisely the set of directions preserving
the condition that the $\psi^e$ are constant, hence is spanned by 
vectors in the form
$C^e = \frac{\partial}{\partial \psi^e_t} - 
\frac{\partial}{\partial \psi^e_s}$. 
So from formula \ref{dif}
we see  that $x$ is  a critical point of the objective 
function in a conformal class
if and only if 
\[0 =  dH(C^e) = \log\left(\frac{\cosh(l^e_t(x))-1}{2}\right) - 
\log\left(\frac{\cosh(l^e_s)-1}{2}\right)\]
for all edges $e$. 
Thus at a critical point of $H$ we have that $l^e_t(x) =l^e_s(x) $, 
so the set of hyperbolic  triangles  formed from
$x$'s angle data
fit together to form an actual
hyperbolic  surface.

\begin{defn}
 A  Delaunay angle system which is the angle data of a hyperbolic surface
 will be called   {\it uniform}.  
\end{defn}

The question becomes: how many (if any) 
uniform structures
can be  associated to a given   
angle system?  The objective function can be analyzed and using  a
 compactness argument with boundary control,  
we find this objective function always
achieves its maximum (see 
\cite{Les}).
$H$  is also observed to be strictly concave down, 
so in fact any critical point is $H$'s unique maximum, and we have

\begin{thm}[The Discrete Uniformization Theorem]
\label{unif}
If $\chi(M) < 0$ and $x$ is a teleportable Delaunay angle system, then $x$ is 
 conformally equivalent to a unique uniform angle system.   
\end{thm}

\begin{remark}
\label{packrem}
At this point it may be unclear how  the above discussion is  
related to a disk pattern problem.
Given a  geodesic triangulation of a hyperbolic surface, 
the disk pattern of interest  here  is the disk pattern produced by 
the circumscribing disks of the
triangulation's  triangles, which we will call the ``empty pattern''. 
Notice in the presence of the ``empty pattern''  we may assign to each edge the
value of the  
intersection angle  
between the circumscribing disks of the triangles sharing this edge.
A pattern production  theorem  in this setting is an 
assurance of the existence of an ``empty pattern''
given a  topological triangulation 
and the specification of a sensible  intersection angle to each edge.
To discover   what sensible means   it is necessary to strengthen lemma
\ref{dtele} to
\begin{lemma}\label{strong}
Given any topological triangulation
and set of data $\psi^e \in
(0,\pi)$ satisfying both  $\sum_{e \in v} \psi^e = 2 \pi$ and 
the teleportability condition, there is $y \in
\frak{N}_x $ satisfying $\psi^e(y) = \psi^e$.
\end{lemma}
Hence from theorem \ref{unif}, under these hypothesis there is 
an ``empty pattern'' with   
intersection angles given by $\pi - \psi^e$, and we have solved our
disk pattern problem.

This   pattern problem is equivalent to a   generalization 
of the convex ideal case of the Thurston-Andreev theorem when
$\chi(M) < 0$ (see \cite{Les} for the details and various
generalizations).  In the Euclidean case, 
this extension was carried out  by Bowditch 
in \cite{Bo}, using techniques similar
to Thurston's original techniques found in \cite{Th}.
The use of an objective function  for solving 
such problems was introduced in 
Colin de Verdi\`ere's \cite{Co} (see Question \ref{det}),
while the  use of hyperbolic volume as an objective function for 
producing disk patterns has its  origin in Br\"agger's beautiful
paper \cite{Be}.  
Hyperbolic volume is also used as an objective function in Rivin's \cite{Ri1}.
\end{remark}

\end{section}

\begin{section}{Continuous Conformal Uniformization}\label{sec:4}

The goal now is to bootstrap from the discrete uniformization
procedure in the previous section to a  procedure  for producing a
conformally equivalent uniform metric on a Riemannian surface.
 Two  metrics on $M$, $g$ and $h$, are conformally 
 equivalent if  $h = e^{2 \phi} g$ for a 
 smooth function $\phi$. 
  For metrics by uniform 
 structure I will mean  a metric with constant
curvature.  Our goal is to prove the classical result...

\begin{thm}[The Metric Uniformization Theorem]
 \label{cuni}
 Every metric is conformally equivalent to a metric 
 of constant curvature, and this metric is unique up to scaling.
  \end{thm}

We will always be thinking in terms of a fixed 
background metric
called $g$ and will label its associated geometric  objects like its gradient,
 Laplacian, curvature, norm, area element, or area  
as $\nabla$,  $\Delta$, $k$,$|\cdot|$ , $dA$ or $A$. 
 For the $h=e^{2 \phi} g$ metric we shall denote these objects with an $h$
subscript.

For starters let us 
note in the metric world we still have  

\begin{lemma}[The Metric Teleportation Lemma]
\label{mtele}
Every metric on a surface 
is conformally equivalent to a metric with  either 
negative, positive, or zero curvature.
\end{lemma}
As in the discrete case, this part of the uniformization procedure 
is completely linear and follows at once from the
facts that $k_{h} = e^{-2\phi}(-\Delta \phi +k)$ and that 
$C^{\infty}(M)$ is the $L^2$ orthogonal direct sum of 
$\Delta(C^{\infty}(M))$
and the constant functions.
With this observation in mind,
in our $\chi(M) <0$ world we will
restrict our attention to metrics with 
strictly negative curvature, and
from here on out we will assume $h$ has
negative curvature everywhere.

As described in section \ref{sec:2}, relative to  a fixed set of vertices 
we may apply Delaunay's ``empty
sphere'' method and produce both the 
Delaunay triangulation and its associated ``empty
disk pattern'' (see remark \ref{packrem}). A conformal 
transformation of a metric preserves 
infinitesimal circles and the angles between them, hence for a dense
enough set of vertices  a conformal transformation of a metric 
nearly preserves the data in the ``empty
disk pattern''. In particular, when we conformally change our metric 
the angles in the associated  
Delaunay triangulation will change in a manner strongly resembling
a discrete conformal change, as introduced in 
the previous section.
With this observation in mind, 
if we choose a dense ${\bf p} = \{p_1, \dots ,p_n\}$ 
and a topological triangulation, $\bf{T}$,   with ${\bf p}$ as its 
vertices, then  we 
could  measure how close to uniform $h$ is  
with the objective function of the previous section. Specifically,
we could let 
\[ {\bf H}_{h}({\bf p}) = \sum_{t \in {\bf T}} V^h_t({\bf p}),\]
where $V^h_t({\bf p})$
is computed using 
the angle data associated to the 
triangulation viewed in the $h$ metric.
This of course  means  connecting the  needed vertices of 
${\bf p}$ with $h$ geodesics 
and measuring the resulting $h$ angles.

To capture this dense enough set of vertices it is
natural to average ${\bf H}_h$ over all sets of vertices 
with a fixed density 
and take the  limit as the vertex
density goes to infinity, as we did for the random variable $F$  in
section \ref{sec:2}.
 To be explicit, we will distribute points with a density $\lambda$ 
relative to  $g$'s area measure,  
and replace the $R_t$ in remark
\ref{unifrem} with the function that is
$V_t^h$ if the ``triple's" disk (in $g$'s metric) is empty of points,
 and zero otherwise.
We will denote the expected value as 
$\Bbb E^g_{\lambda}( {\bf H}_{h})$, with the superscript $g$ there to 
remind us of the
background metric choice.
Let 
\[    I^g(h) = \lim_{\lambda \rightarrow \infty} \Bbb E^g_{\lambda}
({\bf  H}_{h} - {\bf H}_g), \]
where the  second term, ${\bf H}_g$, is independent of $h$ and is needed only 
to normalize the computation. From this construction we 
expect that $I^g$ is an objective function capable of uniformizing 
a negatively curved metric,
and to confirm this,
it is useful to explicitly compute $I^g$.

\begin{thm}\label{form}
If $g$ and $h$ are conformally equivalent with $h = e^{2 \phi}g$ then 
\[ I^g(h)  = -  \int_M | \nabla \phi |^2 + (\Delta \phi - k)
 \log(\Delta \phi - k) + k \log|k| dA.\]
\end{thm}
{\bf Proof:}
In performing this computation,  
we first easily arrive at the analog of 
formula \ref{ex1}, namely  
\begin{eqnarray} 
\label{exn}
I^g(h) =  \frac{1}{6} \int_{V_{\delta}}(V^h(A_1,A_2,A_3)-V^g(A_1,A_2,A_3)) 
e^{-\lambda a(y)} (\lambda dA)^3, 
\end{eqnarray}
with 
$V^h(A_1^h,A_2^h,A_3^h)$  the volume of the prism determined by the
triangle angles $\{A^h_1,A^h_2,A^h_3\}$ formed using the $h$ metric.
At this point theorem \ref{form} follows 
from equation \ref{exn} exactly as equation  \ref{ex3} followed from
equation \ref{ex1}:
We  change  coordinates,
Taylor expand, integrate, and  use the mean value theorem, to arrive
at our formula for $I^g(h)$ plus an error term which in this case is
of 
the form
$O\left(\lambda^{-\frac{1}{2}} \log(\lambda)\right)$.
It should be noted that this procedure, which  
is straightforward when used to  produce equation \ref{ex3}, becomes  
considerably more involved in this setting.
First one must Taylor expand $A^h_i(r,\vec{\theta},p)$ in $r$, which
involves solving the boundary value problem determined by the 
geodesic equation up to the second order.   
With this one may expand $V^h$  in $r$, 
though care is needed  since 
$V^h$'s differential has singularities.  
The singularities can be dealt with, and the 
needed integrals explicitly computed, to 
arrive at theorem \ref{form}. The details can be found in 
\cite{Le}.  
\qed

Now we will use $I^g$
to mimic the 
proof of the discrete uniformization theorem form the previous 
section here in the metric world.

{\bf Proof of Theorem \ref{cuni}}
First we confirm that critical points are uniform by computing 
$I^g$'s differential.  
Using the notation of theorem  \ref{form},
it is natural to view $I^g$ as a function on the
possible $\phi$ where $h = e^{2 \phi} g$.  
With this view point we have
that the Fr\'echet derivative of $I^g$ at $\phi$ in the direction 
$\psi$ is 
\begin{eqnarray} \label{cdif}
D I^g (\psi) 
=   - \int_{M} \Delta \psi \log|k_{h }| dA  .
\end{eqnarray}

{\it Comment.}
Equation \ref{cdif}  implies that the flow 
generated by  $I^g$  is the ``log Ricci" flow, and an alternate proof of 
theorem \ref{cuni} is to use this use this flow as Hamilton used
the Ricci flow in  \cite{Ha}.  


Back to our optimization proof.
A straightforward regularity argument along with equation \ref{cdif} for 
the Fr\'echet derivative
assures us that
$k_{{h}}$ is smooth at a critical point. 
 From formula \ref{cdif} at a critical point
$\log|k_{h}|$ is $L^{2}$ orthogonal to  $\Delta(C^{\infty}(M))$.  
Recalling once again that  $C^{\infty}(M)$ is the $L^2$ 
orthogonal direct sum  of 
$\Delta(C^{\infty}(M))$
and the constant functions, we see that $\log|k_{h}|$ is indeed constant.
So, in
analogy to  the discrete case, a metric is critical if and 
only if it is  uniform.
 
Using a well chosen space of candidate metrics we find that a compactness
argument with boundary control guarantees the existence of a 
critical  point where  $I^g$ achieves its  maximum value.
For example one could use the ``$x \log^{+}(x)$"  Orlicz-Sobolov 
closure (see \cite{Do}) of    
\[V = \{ \phi \in C^{\infty} \mid 
 \int_M \phi dA =0 \mbox{ and } -\Delta \phi + k < 0\}, \]
along with  some basic functional analysis to arrive at 
the needed existence statement (see \cite{Le} for the details).

Since $H_g$ was constructed out of scale invariant 
angle data, we can only hope for uniqueness up to scaling.
As in the discrete case, the  uniqueness follows 
from the concavity of $I^g$.
To be more specific, we use the  fact that all 
critical points are smooth and that the Fr\'echet Hessian at 
$\phi \in V$ applied to 
$ (\psi,\psi)$,
\begin{eqnarray}\label{hess} 
D^2{\bf  }I^g (\psi,\psi) =-  \int_{M} |\nabla \psi|^2  + 
\frac{( \Delta \psi)^2}{\Delta \phi - k} dA,  
\end{eqnarray}
 is strictly negative 
at a non-zero $\psi$ satisfying $\int_M \psi dA =0$.

So we have proved theorem \ref{cuni}  by mimicking the discrete case's 
arguments. 
\qed

Notice that we appear to have an infinite number of objective 
functions, one for each metric $g$.   Fortunately any pair of 
these objective functions differ only 
by a constant. In order to see this it is useful to recall two other 
functions related to 
metric uniformization, the  $\log(\det(\Delta_h))$ and the metric
entropy. 
The $\log$ of the determinant of the Laplacian 
 had its  uniformization properties explored by 
Osgood, Phillips, and Sarnak in \cite{os}, while the entropy
\[ E(h) = -\int_M k_h \log|h_h| dA_h \]
turned up Hamilton's  paper   on surface 
 uniformization  \cite{Ha}.
As a straightforward consequence of theorem  \ref{form}
 and  Polyakov's formula (see \cite{Po})   telling us that  the
$\log(\det(\Delta_h))$ with in a conformal class can be expressed as 
\[ \log(\Delta_h) = - \frac{1}{6 \pi} 
\left(\frac{1}{2} \int_M ||\nabla \phi||^2 dA + \int_{M} k \phi dA  
+ \ln(A_h) \right) +C(g), \]
we have

\begin{cor}\label{newform}
Let $g$ and $h$ be conformally equivalent and let 
 \[M(h) = E(h)-
12 \pi \log\left(\frac{\det(\Delta_h)}{A_h}\right).\] 
Then $ I^g(h) = M(h) - M(g)$.
\end{cor}

So up to a constant our objective function 
is given by $M(h)$.
\end{section}

\begin{section}{Questions}\label{sec:5}

\begin{question}\label{det}
In \cite{Co}, Colin de
Verdi\`ere suggested that objective functions
related to circle pattern problems
might be related to  the determinant of the Laplacian.
One certainly would hope for a more direct relationship
than that provided by corollary \ref{newform}, raising the following natural
question:  Is there a disk pattern objective function,
${\bf D}_h$,  with the property that
$\Bbb E^g_{\lambda} ({\bf  D}_{h})$ limits to 
$\log(\det(\Delta_h))$ when suitably normalized?
\end{question}

\begin{question}\label{nonconformal}
Corollary \ref{newform}  also reveals some unexpected properties of $I^g$.
For example, there was no reason to expect that the
roles of the $g$ and $h$ metric could be decoupled via $I^g(h) = 
M(h)-M(g)$; or  even the
mysterious implication that 
$I^g(h)=-I^h(g)$.  
Furthermore corollary \ref{newform} allows
us to extend  this objective function consistently to the space of all
negatively curved metrics. 
Notice that $I^g(h)$ can be computed and discretely interpreted 
when $g$ and $h$ are not in the same
conformal class.  With this observation we are left with the 
following natural question: 
Does  corollary \ref{newform}  hold  among all metrics with
negative curvature?  Notice in particular this would provide a local 
formula for the 
$\log(\det(\Delta_h))$ outside a conformal class.
\end{question}

\begin{question}
The idea of using disk patterns to explore  uniformization can be
traced back to Thurston  \cite{Thu}.  Thurston's original
idea was to approximate the Riemann mapping with a disk pattern
solution.   Thurston's idea was initially justified in \cite{Su} and
has been developed considerably  since then, see for example \cite{He}.
It would be interesting to implement the spirit of this approximation
approach in this setting. For example it would be nice
to answer the following question: If one takes a ``dense'' set of
points on a Riemannian surface and  forms the "empty disk pattern",
can one  measure in a meaningful way how far the uniform surface
produced by theorem \ref{unif} is from the conformally equivalent
uniform structure produced by theorem \ref{cuni}? 
\end{question}

\begin{question}
This whole story carries over to the spherical case, with the exception
of the convexity of the analogs of the objective functions $H$ and $I^g$.
Can one use these same techniques to prove the corresponding
uniformization results, even without convexity?
\end{question}

\begin{question}\label{eyesontheprize}
Can techniques like these be applied to find geometric structures
on 3-manifolds?
\end{question}

\end{section}


\begin{thebibliography}{10}

\bibitem{Bo}
B.~H. Bowditch.
\newblock Singular euclidean structures on surfaces.
\newblock {\em J. London Math. Soc. (2)}, 44:553--565, 1991.

\bibitem{Be}
W.~Br{\"a}gger.
\newblock Kreispackungen und triangulation.
\newblock {\em Ens. Math}, 38:201--217, 1992.

\bibitem{Co}
Y.~Colin de~Verdi\'ere.
\newblock Un principe variationnel pour les empilements de cercles.
\newblock {\em Pr\'epublication de l'Institut Fourier, Grenoble}, 147:0--16,
  1990.

\bibitem{Da}
B.~Delaunay [Delone].
\newblock Sur la sphe\'re vide.
\newblock {\em Proc. Internal. Congr. Mth.}, 1:695--700, 1928.

\bibitem{Do}
T.~Donaldson.
\newblock Non-linear elliptic boundary problems in orlicz--sobolov spaces.
\newblock {\em J. Diff. Eq.}, 10:507--528, 1971.

\bibitem{Ha}
R.S. Hamilton.
\newblock The ricci flow on surfaces.
\newblock {\em Contemporary Math.}, 71:237--262, 1988.

\bibitem{He}
Z.~He and O.~Schramm.
\newblock The $c^{\infty}$-convergence of hexagonal disk packings to the
  riemann map.
\newblock {\em Acta Math.}, 180:219--245, 1998.

\bibitem{Les}
G.~Leibon.
\newblock Singular hyperbolic structures on a surface.
\newblock In Preparation.

\bibitem{Le}
Gregory Leibon.
\newblock {\em Random Delaunay Triangulations, the Thurston-Andreev Theorem and
  Metric Uniformization}.
\newblock PhD thesis, UCSD, 1999.

\bibitem{Po}
A.~Polyakov.
\newblock Quantum geometry of bosonic strings.
\newblock {\em Phys. Lett. B}, 103:207--210, 1981.

\bibitem{os}
B.~Osgood R.~Phillips and P.Sarnak.
\newblock Extremals of determinants of laplacians.
\newblock {\em J. Fun. Anal}, 80:148--211, 1988.

\bibitem{Ri1}
I.~Rivin.
\newblock Euclidean structures on simplicial surfaces and hyperbolic volume.
\newblock {\em Ann. of Math (2)}, 139:553--580, 1994.

\bibitem{Su}
B.~Rodin and D.~Sullivan.
\newblock The convergence of circle packings to the riemann mapping.
\newblock {\em J. Differential Geom.}, 26:349--360, 1987.

\bibitem{Th}
William~P. Thurston.
\newblock {\em Three-Dimensional Geometry and Topology}.
\newblock The Geometry Center, University of Minnesota, draft edition, 1991.

\bibitem{Thu}
W.P. Thurston.
\newblock The finite riemann mapping theorem.
\newblock Unpublished talk given at the International Symposium in Celebration
  of the proof of the Bieberbach Conjecture (Purdue University, 1985).

\end{thebibliography}
\end{document}